\documentstyle[diagram,12pt]{article}

\newcommand{\ts}{\otimes}

\newcommand{\qq}{{\bf Q}}
\newcommand{\zz}{{\bf Z}}

\newcommand{\oo}{{\cal O}}
\newcommand{\ii}{{\cal I}}
\newcommand{\jj}{{\cal J}}

\newcommand{\ra}{\rightarrow}

\newcommand{\ord}{{\rm ord}}
\newcommand{\red}{{\rm red}}
\newcommand{\base}{{\rm BS}}
\newcommand{\mult}{{\rm mult}}
\newcommand{\supp}{{\rm supp}}
\newcommand{\codim}{{\rm codim}}

\newtheorem{theorem}{Theorem}[section]
\newtheorem{lemma}[theorem]{Lemma}

\newtheorem{conjecture}[theorem]{Conjecture}
\newtheorem{corollary}[theorem]{Corollary}
\newtheorem{proposition}[theorem]{Proposition}
\newtheorem{definition}[theorem]{Definition}
\newcounter{nmb}

\newcommand{\next}{\addtocounter{nmb}{1}}
\newcommand{\rn}{\addtocounter{theorem}{1}}
\newcommand{\nmbthm}{\renewcommand{\theequation}{\thetheorem.\thenmb}}
\newcommand{\nonmbthm}{\renewcommand{\theequation}{\thetheorem}}
\newcommand{\init}{\setcounter{nmb}{0}}
\newenvironment{thm}{\setcounter{nmb}{0}\begin{theorem}}{\end{theorem}}
\newenvironment{lm}{\setcounter{nmb}{0} \begin{lemma}}{\end{lemma}}

\newenvironment{define}{\setcounter{nmb}{0}\begin{definition}}
{\end{definition}}

\begin{document}
\setcounter{section}{-1}
\title{Base loci of linear series are numerically determined}
\author{Michael Nakamaye}
\maketitle

\section{Introduction}
\nonmbthm

Suppose $X$ is a smooth projective variety defined over an algebraically
closed field of characteristic zero and suppose $L$ is a line bundle on $X$.
If 
$x \in X$ the Seshadri constant 
$\epsilon(x,L)$ measures the numerical positivity of $L$ at $x$.  
When $L$ is nef $\epsilon(x,L)$ carries local geometric information at $x$,
namely it measures which order jets $L$ generates at $x$ asymptotically.
When $L$ is not nef, however, $\epsilon(x,L)$ often
ceases to carry geometrically
meaningful information.  In order to remedy this situation
we study a slight variant of the Seshadri constant, 
$\epsilon_m(x,L)$, which we call the moving Seshadri constant of $L$ at $x$
and which specifically  measures which order
jets a high  multiple of $L$ generates at $x$.  The main result of this 
paper is to generalize the results of \cite{N}, establishing that the base
locus of a small perturbation of a big line bundle $L$ is completely determined
by the constants $\epsilon_m(x,L)$.

\medskip

We begin by recalling the basic 
definitions and results of \cite{N}.

\begin{define}  
Suppose $L$ is a line bundle on a smooth projective variety $X$.
Let
$$
\base(L) = \{x \in X: \mbox{$s(x) = 0$ for all sections
$s \in H^0\left(X,nL\right)$  and for all $n > 0$}\}.
$$
We call $\base(L)$ the stable base locus of $L$.
\label{def1}
\end{define}
Note that 
$\base(L)$ is empty  if and only if some positive multiple of $L$ is
generated by global sections.  
If $D \in {\rm Div}(X) \ts \qq$ is a divisor with rational coefficients, then
Definition \ref{def1} can be slightly reformulated:
$$
\base(L(-D)) = \base\left(kL(-kD)\right)
$$
where $kD \in {\rm Div}(X)$.  This definition is
independent of $k$.  
\begin{define}
Suppose $L$ is a nef line bundle on $X$.  Let
$$
L^\perp = \left\{ P \in X: c_1(L)^{\dim V}\cap V = 0 \,\,\,\mbox{for some}\,\,V
\subset X\,\,\mbox{containing}\,\,P\right\}.
$$ 
\label{def2}
\end{define}
We established in \cite{N} that $L^\perp$ is Zariski closed. 
The main theorem of \cite{N} is the following:
\begin{thm}  Suppose $X$ is a smooth projective variety of dimension $d
\geq 2$ and $L$ a big and nef line bundle on $X$.  Suppose $A$ is an ample
line bundle on $X$.  Then there exists $\epsilon > 0$ such that 
$$
L^\perp = \base(L(-\delta A))
$$
whenever $0 < \delta < \epsilon$.
\label{main}
\end{thm}

Thus the base locus of a nearly big and nef line bundle is determined 
numerically.  Here we push this analysis further and show that at each stage
where the base locus jumps, there is a numerical reason for the jump which 
can be captured by studying moving Seshadri constants which we define here:

\begin{definition}
Suppose $X$ is a smooth projective variety and $L$ a big line bundle on $X$.
For all $n \gg 0$, let $\pi_n: X_n \ra X$ be a resolution of $\ii_{B_n}$ where
$B_n$ is the base locus scheme of the complete linear series $|nL|$. We
assume here that $X_n$ is smooth.  Write
$$
\pi_n^\ast (nL) \sim M_n + E_n
$$
where $\oo_{X_n}(-E_n) = \pi_n^{-1}\left(\ii_{B_n}\right)\cdot \oo_{X_n}$.
For any point $x$ not in $\base(L)$, we define the moving Seshadri constant
of $L$ at $x$ by
$$
\epsilon_m(x,L) = \lim_{n \ra \infty}\frac{\epsilon(x,M_n)}{n}.
$$
\label{mvsesh}
\end{definition}
The invariant $\epsilon_m(x,L)$
measures how many leading term can be specified at $x$, 
asymptotically, in the Taylor
series of a section $s \in H^0(X,nL)$.  
Note that the limit in Definition \ref{mvsesh} makes sense.
Indeed, since 
$$
\epsilon(x,M_n) \leq \sqrt[\frac{1}{\dim X}]
{c_1(M_n)^{\dim X}},\,\,\,\,\forall
n,
$$
and
$$
\frac{c_1(M_n)^{\dim X}}{n^{\dim X}\dim X!} \leq
\lim_{m \ra \infty}\frac{h^0(X,mnL)}{(mn)^{\dim X}}
$$
it follows that 
$a = \limsup_{n \ra \infty}\left\{
\frac{\epsilon(x,M_n)}{n}\right\}$ exists.
Suppose there were an infinite sequence $\{n_1,n_2,\ldots\}$ such that 
\rn
\begin{eqnarray}
\lim_{m \ra \infty} \frac{\epsilon(x,M_{n_m})}{n_m} = b < a.
\label{yy}
\end{eqnarray}
By \ref{yy}, we can choose a 
positiive integer $k$ such that 
\rn
\begin{eqnarray}
\epsilon(x,M_k) > b.
\label{zz}
\end{eqnarray}
For any $r$ 
we can write
$n_r = ck + d,\,\, d < k$. 
Using the fact that
$$
\frac{\epsilon(x,M_{ck+d})}{ck+d} \geq \frac{\epsilon(x,cM_k)+
\epsilon(x,M_d)}{ck+d}
$$
we find that 
$$
\epsilon(x,M_{n_r}) \geq  \frac{c\epsilon(x,M_k)+
\epsilon(x,M_d)}{ck+d}
$$
Since $d$ is bounded and $c \ra \infty$ as $r \ra \infty$ this contradicts
\ref{zz}. 

Next, we isolate those line bundles where the stable base locus changes.
\begin{definition}
Suppose $D$ is a ${\bf Q}$--divisor on $X$ in the interior of the
effective cone, $A$ an ample $\qq$--divisor.  Suppose
$\alpha$ is a positive real number.  We call $\alpha$ a jumping number relative
to $D$ and $A$ if 
$\base(D-\alpha A+\epsilon_1 A)$ is properly contained in
$\base(D- \alpha A -\epsilon_2 A)$
for any $\epsilon_1, \epsilon_2 > 0$ such that $\alpha +\epsilon_1 $ and
$\alpha -\epsilon_2 $ are rational numbers.
\end{definition}
Intuitively, the stable base locus of $D-\beta A$ changes each time
$\beta$ passes
through a jumping value $\alpha$.  
Our main theorem is to show that jumping numbers are determined by
the relevant moving Seshadri constants approaching zero:

\begin{theorem}
Suppose $E$ is a big $\qq$--divisor on a smooth projective variety $X$ and $A$
an ample divisor.  Suppose $\alpha$ is a jumping number relative to $E$ and
$A$.  We will write
$$
D = E - \alpha A.
$$
Let $Y$ be an irreducible component of $\base(D-\beta
A)$ for $0 < \beta < \epsilon$  
such that $Y$ is not in  $\base (D + \beta  A)$ for any $\beta > 0$.  
 Then for a general  point $\eta \in Y$
$$
\lim_{\beta \ra 0 } \epsilon_m(\eta,D + \beta A) = 0;
$$
here the limit is of course from above as $\epsilon_m(\eta,D + \beta A)$ 
has not been defined for $\beta \leq 0$.  
\label{seshlimit}
\end{theorem}

Note the similarity between Thereom \ref{seshlimit} and Theorem \ref{main}.
In both cases, the critical stage where the base locus changes is marked
by the appropriate Seshadri constants approaching zero.  The proof of
Theorem \ref{seshlimit} is  more subtle than that of Theorem \ref{main}
because the divisors $D + \beta A$, as $\beta$ approaches 
$0$, are not nef and their base loci cannot necessarily
all be resolved on a single birational model of $X$.
Moreover, and more seriously, the critical number $\alpha$ may well be
irrational making it impossible to actually consider a limit linear series. 
The proof of Theorem \ref{seshlimit} nonetheless
follows along the same lines as the 
proof of Theorem \ref{main}, namely by lifting a non--vanishing
section on a subscheme supported on $Y$ using a 
cohomological vanishing theorem.  
In order to avoid working on infinitely many birational models of $X$
 we examine and compare
carefully the base locus scheme of $n(D +\beta A)$
for $\beta$ slightly less than $0$ and for $\beta$ slightly greater
than $0$.  Theorem \ref{seshlimit} amounts to a type of ``continuity''
of these base loci.

\medskip

The outline of this paper is as follows.  In \S 1 we prove Theorem 
\ref{seshlimit}.  
The ideas employed in the proof are very close to the work of Demailly, Ein,
and Lazarsfeld \cite{DEL} on asymptotic multiplier ideals 
but the key technical 
lemmas require other limiting procedures and we were unable
to obtain these results mechanically from \cite{DEL}.
The fundamental idea is simply to keep track of how the base loci of
$|nD|$ vary under small perturbations and how cohomology groups on these
base loci vary.  The key insight is to restrict ourselves to a very specific
vector space of sections, namely those generating specific jets at a general
point $\eta \in Y$.  
In \S 2 we will apply the techniques developed in \S 1 to prove a theorem
of Fujita, also obtained in \cite{DEL}, in a very concrete fashion which
should be helpful in many applications.  

\medskip

Our original goal in this work was to obtain a sufficiently deep understanding
of stable base loci and the effective cone in order 
to derive certain consequences in diophantine
geometry.  The present work comes close but needs to be made more effective
before yielding the desired corrolaries.  Nevertheless, 
the problem of understanding the
effective cone, and defining numerical invariants on it, seems sufficiently
interesting in its own right to be pursued.  In particular, if one defines
$\epsilon(x,L)$ to be zero whenever $x \in \base(L)$ then $\epsilon(x,L)$
can be viewed as a function on the rational classes in the
cone of divisors modulo numerical equivalence:  in this context, the
content of Theorem \ref{seshlimit} is to establish that this function
is continuous with respect to $L$.  

\medskip
\noindent {\em Notation and Conventions}
\begin{itemize}
\item  If ${\cal F}$ is a coherent sheaf on $X$ then $h^i(X,{\cal F}) =
{\rm dim}\,H^i(X,{\cal F})$.
\item If $\ii \subset \oo_X$ is an ideal sheaf
then we let $Z(\ii) \subset X$ denote the associated subscheme of zeroes of
$\ii$.  If $f$ is a rational function on $X$ we denote by $Z(f)$ the
divisor of zeroes of $f$.  
\item  If $D$ is a Cartier divisor on $X$ we will often abbreviate
$H^0(X,\oo_X(D))$ by $H^0(X,D)$.
\item  If $V \subset X$ is a subvariety we denote by $\ii_V^a$ the ideal
sheaf of functions whose order of vanishing along $V$ is at least $a$.  
Similarly we will write $aV$ for the scheme defined by the ideal sheaf
$\ii_V^a$.  
\item If $Y$ is a scheme we denote by $Y_{\rm red}$ the reduced scheme
with support equal to the support of $Y$.
\item  For a variety $X$ we let $K(X)$ denote the field of rational 
functions on $X$.
\end{itemize}

\section{Continuity of Seshadri Constants}

Choose an $\epsilon > 0$ so that $Y$ is an irreducible component
of $\base(D -\beta A)$ for all $0 < \beta \leq \epsilon$.  
We also assume that $D - \epsilon A \in {\rm Div(X)} \ts \qq$.  
Next
choose an $n >0$ so that $n(D-\epsilon A)$ is a genuine line bundle
whose sections give a birational map to projective space. Let
$$
\pi: X^\prime \ra X
$$
be a birational map, with $X^\prime$ smooth, such that
$$
\pi^\ast (n(D-\epsilon A)) \sim M_n + E_n
$$
where $M_n$ is generated by global sections and $\pi^{-1} \ii_B \cdot 
\oo_{X^\prime} = \oo_{X^\prime}(-E_n)$ with $B$ denoting the base locus 
scheme of $|n(D-\epsilon A)|$.  Since $M_n$ is big and nef we
have 
\rn
\begin{eqnarray}
h^1(X^\prime,kM_n) \leq O\left(k^{\dim X -1}\right).
\label{first}
\end{eqnarray}

Let $Z_n$ be the zero scheme of $\pi_\ast \oo_{X^\prime}(-E_n)$ and let
$kZ_n$ denote  the zero scheme of the ideal sheaf
$$
\pi_\ast\oo_{X^\prime}(-kE_n) \subset \oo_X.
$$  
Note that
by \cite{ev1} 3.3, $Z_n$ does not depend on the choice of resolution $\pi$.
Using the long exact cohomology sequence associated to the exact sequence
of sheaves
\begin{eqnarray*}
0 \ra \ii_{kZ_n} \ts \oo_X(kn(D-\epsilon A)) &\ra&
\oo_X(kn(D-\epsilon A))\\
&\ra& \oo_{kZ_n}(kn(D-\epsilon A)) \ra 0
\end{eqnarray*}
gives
\rn
\begin{eqnarray}
H^0(X,\oo_X(kn(D-\epsilon A))) &\ra& H^0(kZ_n,
\oo_X(kn(D-\epsilon A))) \nonumber \\
&\ra& H^1(X,\oo_X(kn(D-\epsilon A) \ts \ii_{kZ_n})) .
\label{thirdy}
\end{eqnarray}
By \ref{first}, pushed down to $X$
$$
h^1(X,\oo_X(kn(D-\epsilon A) \ts \ii_{kZ_n})) = O\left(k^{\dim X
-1}\right)
$$
and thus most sections of $H^0(kZ_n,
\oo_X(kn(D -\epsilon A))) $ lift to $X_n$.  In order to prove Theorem
\ref{seshlimit} we will show that if 
$$
\lim_{\beta \ra 0 } \epsilon_m(\eta,D+ \beta A) \neq 0
$$
then it is possible to use \ref{thirdy} to lift to $X$ a section of 
$H^0(kZ_n,\oo_X(kn(D-\epsilon A))) $ with small order of vanishing
along $Y$.  On the other hand, we will establich that all sections of
$H^0\left(X,\oo_X(kn(D-\epsilon A))\right)$ 
have large order of vanishing along $Y$, thus deriving a contradiction.  

\medskip
Since $\epsilon_m(\eta,D +\beta A)$ is a decreasing function of $\beta$,
if Theorem \ref{seshlimit} were false, there must exist a $\delta > 0$ such
that 
\rn
\begin{eqnarray}
\epsilon_m(\eta,D +\beta A) > \delta ,\,\,\,\, \forall \beta > 0.
\label{fifth}
\end{eqnarray}
We will use the following four results in the proof of Theorem \ref{seshlimit}.
\begin{lemma} With notation as above, there exists $\gamma > 0$, not depending
on $\epsilon$ or $n$, so that
$$
\ii_{Z_n} \subset \ii_Y^{\gamma \epsilon n} ,\,\,\,\,\forall n \,\,\,\mbox{
sufficiently divisible}
$$
and, in fact, $\ii_{kZ_{n}} \subset 
\ii_Y^{\gamma \epsilon kn}$ for all $k>0$.
\label{l1}
\end{lemma}

Lemma \ref{l1} simply states that once a subvariety $Y$ enters the stable
base locus of $D$,
multiplicity along $Y$ must grow linearly as $D$ moves away from the ample
cone.  This result is a formal consequence of properties of multiplier ideal
sheaves when phrased in the language in \cite{DEL}; we give an alternative
proof using \cite{eln}.

Write
\rn
\begin{eqnarray}
Z_n = Y_n \cup W_n
\label{sixth}
\end{eqnarray}
where $Y_n$ is supported on $Y$ and $\supp(W_n)$ does not contain $Y$.
To avoid ambiguity, one can choose $\ii_{W_n}$ maximal over all such
possible expressions; in other words, the local primary
decomposition of $\ii_{W_n}$ should have no embedded components supported
on $\supp(Y)$.
In order to prove Theorem \ref{seshlimit}, we will need to consider the ideal
sheaves $\ii_{Z_n}$ as depending on $\epsilon$ and then take the limit as 
$\epsilon \ra 0$.  We will write $Z_n(\epsilon)$, 
$Y_n(\epsilon)$, and
$W_n(\epsilon)$ in this case, including the possibility that $\epsilon < 0$;
note that $n$ actually depends upon
$\epsilon$ as it is the multiple used to clear denominators in $\qq$ divisors.

Let $\phi: V \ra X$ be a resolution of the base
locus of $|n(D+\epsilon A)|$ with exceptional divisor $E$ so
that $\ii_{Z_n( - \epsilon)} = \phi_\ast \oo_{V}(-E)$ by definition.

\begin{lemma} Assume that $\gamma \epsilon < \delta$.
For all $k \gg 0$ and for all $n$ sufficiently large
there is a vector subspace
$$
W_{k,n}(\epsilon) \subset
H^0\left(kY_{n}(\epsilon),
kn(D + \epsilon A)\ts \ii_{kZ_n(-\epsilon)}\right)
$$
with
$$
\dim(W_{k,n}(\epsilon)) \geq
O\left((\gamma \epsilon nk)^{\codim (Y,X)}(\delta nk)^{\dim Y}\right)
$$
and such that no non--zero section $s \in W_{k,n}$ vanishes
 to order $\geq \gamma \epsilon
kn$ along $Y$.
Moreover, the implied constant depends only on $X$.
\label{lemma2}
\end{lemma}

Lemma \ref{lemma2} merely states that there will exist a large space of
functions on $Y$ which vanish along $kZ_n(- \epsilon)$; this is
clear since we have a lower bound on the Seshadri constant of $D +
\epsilon A$ 
along a generic point of $Y$, independent of $\epsilon$,
 and so this space of sections will
come from restricting sections of $H^0(X,kn(D + \epsilon A))$ which generate 
appropriate jets at $\eta$.  
For the following Corollary, we assume for simplicity that $D \in 
{\rm Div}(X) \ts \qq$ and hence that $\epsilon \in \qq$.  We will show how to
avoid this assumption in the proof of Theorem \ref{seshlimit}.

\begin{corollary}
With the same notation as in Lemma \ref{lemma2}, let $B$ be an ample line
bundle on $X$ and for $a \in (0,1) \cap \qq$, let
 $V_{k,n}(a\epsilon) \subset W_{k,n}(a\epsilon)$
be defined by
$$
V_{k,n}(a\epsilon) = 
\{s \in W_{k,n}(a\epsilon)| Z(s)- akn B \,\,\,\, \mbox{is effective}
\}.
$$
Then there exists an integer $n$ such
that for all $k\gg0$
$$
\dim(V_{k,n}(a\epsilon)) \geq
O\left((\gamma a\epsilon nk)^{\codim (Y,X)}(\delta nk)^{\dim Y}\right) - 
O\left(a^{\codim (Y,X)+1} (nk)^{\dim X}\right)
$$
with implied constants depending only on $X$, $\epsilon$, and $B$.  
\label{cor}
\end{corollary}

\noindent
Corollary \ref{cor} examines what happens to the subspace of sections
of Lemma \ref{lemma2} when we replace the divisor
$kn(D+ \epsilon A)$ with the less positive divisor
$kn(D- \epsilon A)$, the idea being to limit how many of the
sections of Lemma \ref{lemma2} can be lost.

For our final Lemma, we adopt the following notational conventions.  We let
$n$ denote an integer to be fixed in the proof.  For each positive integer
$m$ we let 
\begin{eqnarray*}
n_m = 2^mr_mn, \\ 
n_m^\prime = (2^m-1)r_mn,
\end{eqnarray*}
where $r_m$ is a positive integer to be specified in the proof. 
\begin{lemma} Fix a  number $0 < \beta < \epsilon$ so that $D - \beta A \in
{\rm Div(X)} \ts \qq$.
There exists an ideal sheaf ${\cal J}$ 
on $X$, with $Y$ not contained in $Z(\jj)$,
 such that for all $k$ sufficiently large
and all $m \geq 0$
$$
\jj^{\frac{kn_m}{2^m}} \cdot \ii_{kZ_{n_m^\prime}(- \beta)} \subset 
\ii_{kW_{n_m}(- \beta + \frac{\epsilon}{2^m})}
$$
Moreover, if $D \in {\rm Div(X)} \ts \qq$ 
then this relation also holds for $\beta = 0$. 
\label{lemma3}
\end{lemma}

Having moved from 
$kn(D + \epsilon A)$ to
$kn(D -\epsilon A)$ with Corollary \ref{cor} we now use Lemma
\ref{lemma3} to replace
$\ii_{kZ_n(- \epsilon)}$ with 
$\ii_{kW_n(\epsilon)}$ in Lemma \ref{lemma2}.  
Combined with the previous results and
the exact sequence \ref{thirdy},
this will exhibit a large vector space of sections of
$H^0(X,kn(D -\epsilon A))$ having small order of vanishing
along $Y$, contradicting Lemma \ref{l1}.
More specifically,

\medskip

\noindent
{\bf Proof of Theorem \ref{seshlimit}}  
\nmbthm
\setcounter{theorem}{9}
We will assume, for simplicity, that $D \in {\rm Div}(X) \ts \qq$ 
and show how to remove this hypothesis at the end.  
Applying Lemma 
\ref{lemma3} for $\beta = \frac{\epsilon}{2^{m}}$, and replacing $\jj$ by
$\jj^2$, 
gives an ideal sheaf
$\jj$ satisfying, for $k$ sufficiently divisible relative to $m$
\next
\begin{eqnarray}
\jj^{\frac{kn_m}{2^m}} \cdot \ii_{kZ_{n_m^\prime}(- 
\frac{\epsilon}{2^{m}})}
 \subset 
\ii_{kW_{n_m}(\frac{\epsilon}{2^{m}})} ,\,\,\,\,\forall m \geq 0.
\label{1st}
\end{eqnarray}
Choose a line bundle $C$ on $X$ so that $C \ts \jj$ is generated by global 
sections.  We let
$$
B = C + 2\epsilon A.
$$
Applying Corollary \ref{cor} with $a = \frac{1}{2^m}$ gives a subspace
$$
U_{k,n_m^\prime}(a\epsilon) \subset
H^0\left(kY_{n_m^\prime}(a\epsilon), kn_m^\prime \left(D + a\epsilon A \right)
\ts \ii_{kZ_{n_m^\prime}}(-a\epsilon)\right)
$$
with bounded dimension:
\next
\begin{eqnarray}
\dim(U_{k,n_m^\prime})(a\epsilon) \geq
O\left((\gamma a\epsilon n_m^\prime k)^{\codim (Y,X)}
(\delta n_m^\prime k)^{\dim Y}\right) - 
O\left(a^{\codim (Y,X)+1} (nk)^{\dim X}\right).
\label{bb2}
\end{eqnarray}
Since all sections in $U_{k,n_m^\prime}(a\epsilon)$
are defined via restriction from sections on $X$ 
and since none of these sections vanishes on $\gamma a \epsilon
kn_m^\prime Y$
we can consider
$$
U_{k,n_m^\prime}(a\epsilon) \subset
H^0\left(kY_{n_m}(a\epsilon), kn_m^\prime \left(D + a\epsilon A \right)
\ts \ii_{kZ_{n_m^\prime}}(-a\epsilon)\right).
$$

By hypothesis, every section $s \in U_{k,n_m^\prime}(a\epsilon) $
vanishes along $B = C + 2 \epsilon A $ 
to multiplicity at least $akn_m^\prime$; removing the base divisor $B$
gives a vector space
$$    
V_{k,n_m^\prime}(a\epsilon) \subset
H^0\left(kY_{n_m}(a\epsilon), kn_m^\prime \left(D - a \epsilon A - a C \right)
\ts \ii_{kZ_{n_m^\prime}}(-a\epsilon)\right).
$$
Note that removing the base divisor $B$ does not
effect membership in the ideal sheaf $\ii_{kZ_{n_m^\prime}}(-a\epsilon)$. 
Tensoring the vector space $V_{k,n_m^\prime}(a\epsilon)$ 
by $akn_m^\prime C$
and using \ref{1st}, taking care again to replace 
$\jj$ with $\jj^2$ to account for the difference between $n_m$ and 
$n_m^\prime$,
 allows us to view $V_{k,n_m^\prime}(a\epsilon)$ as
$$
V_{k,n_m^\prime}(a\epsilon) \subset
H^0\left(kY_{n_m}(a\epsilon), kn_m^\prime \left(D - a \epsilon A \right)
\ts \ii_{kW_{n_m}}(a\epsilon)\right).
$$
Note that since $\eta \not\in \supp(C)$
no section of $V_{k,n_m^\prime}(a\epsilon)$ 
has multiplicity greater than 
$\gamma \epsilon akn_m^\prime$ at $\eta$.  
We now replace $\gamma$ with $\gamma/2$
and choose a section 
$$
\sigma \in H^0\left(X, k(n_m - n_m^\prime) (D - a \epsilon A))\right)
$$
of minimal index at $\eta$ and tensor to obtain an inclusion
\next
\begin{eqnarray}
 V_{k,n_m^\prime}(a\epsilon) \stackrel{\ts \sigma}{\hookrightarrow}
H^0\left(kY_{n_m}(a\epsilon), kn_m \left(D - a \epsilon A \right)
\ts \ii_{kW_{n_m}}(a\epsilon)\right).
\label{bb4}
\end{eqnarray}
For $m$ sufficiently large, 
every section in the vector subspace \ref{bb4} has multiplicity at most
$ \gamma a\epsilon kn_m$ at $\eta$: here we use the extra $\gamma /2$ 
obtained by replacing $\gamma$ by $\gamma /2$.

Using \ref{bb2} we see that for $m$ sufficiently large
$$
\dim (V_{k,n_m^\prime}) \geq
O\left((\gamma a\epsilon nk)^{\codim (Y,X)}(\delta nk)^{\dim Y}\right) 
$$
But by \ref{sixth} 
any non--zero section $s$  in the vector subspace \ref{bb4} glues together
with the zero section on $kW_{n_m}(a\epsilon)$
to give a non--zero section of
$$
H^0\left(
kZ_{n_m}(a\epsilon),
kn_m\left(D -a\epsilon A\right)\right).
$$
Using \ref{thirdy}, however, we obtain, for $m$ sufficiently large,
many non--zero sections of 
$$
H^0\left(
X,kn_m\left(D -a\epsilon A\right)\right)
$$
with order of vanishing strictly less than
$\gamma a\epsilon kn_m$ along $Y$.  This contradicts Lemma \ref{l1},
or more precisely \ref{order},
and establishes Theorem \ref{seshlimit} in case $D \in {\rm Div}(X) \ts \qq$.

\medskip
Suppose now that $D$ is only a real divisor and consequently that
$\epsilon$ is irrational; in particular, Corollary \ref{cor} does not hold
as stated and  Lemma \ref{lemma3} does not
hold for $\beta = 0$.
The same proof holds except that we need to employ an additional
limiting procedure.
First, for \ref{1st}
which no longer makes any sense, we choose 
$$
\beta_m, \gamma_m \in \left(\frac{\epsilon}{2^{m+1}},\frac{\epsilon}
{2^m}\right)
$$
so that $D - \beta_m A, D + \gamma_m A \in {\rm Div}(X) \ts \qq$.  
Then for $k$
sufficiently divisible with respect to $m$ we have
$$
\jj^{\frac{kn}{2^m}} \cdot \ii_{kZ_{n_m^\prime}(- \beta_m)} \subset 
\ii_{kW_{n_m}(\gamma_m)} ,\,\,\,\,\forall m \geq 0.
$$
As far as Corollary \ref{cor} goes, the problem is that $W_{k,n_m}(a\epsilon)$
will not make sense in general because $D + a\epsilon A$ may not be a rational
class.  This is not a problem, however, since for any rational class
$D + b A$ with $b > a\epsilon$, the estimate of
Corollary \ref{cor} will still hold so the same limiting
procedure applies.

\medskip
\noindent
{\bf Proof of Lemma \ref{l1}}
\nonmbthm
\setcounter{theorem}{4}
\init
\nmbthm
Suppose
$\phi: V \ra X$ is a birational map with $V$ smooth and suppose that there
is a unique smooth exceptional divisor
$E_Y$ dominating $Y$: this can be obtained, for example, by taking an embedded
resolution of $Y$ in $X$ followed by a blow--up of the strict transform of
$Y$.
Consider the commutative diagram
$$
      \begin{diagram}
	\node{V^\prime}\arrow[1]{e,t}{\phi^\prime}
	\arrow{s,l}{\pi^\prime}\node[1]{X^\prime}
	\arrow{s,r}{\pi}\\
	\node{V}\arrow[1]{e,t}{\phi}
	\node[1]{X}
      \end{diagram}
$$
Setting
$E_n^\prime = (\phi^\prime)^\ast E_n$, we 
claim that it is sufficient to show that
\next 
\begin{eqnarray}
E_n^\prime - (\pi^\prime)^\ast
\gamma \epsilon n \tilde{E}_Y \,\,\,\,\mbox{is effective}.
\label{t1}
\end{eqnarray}
To see this suppose that $\eta$ is a general point of $Y$.
If for some $k>0$, $\ii_{kZ_{n}}$ were not contained in 
$\ii_Y^{\gamma \epsilon kn}$ then
there would exist an open subset $U \subset X$ 
and a function $f \in \Gamma(U,\oo_U)$
such that $\ord_\eta(f) < \gamma \epsilon kn$ but $(\phi^\prime)^\ast
(\pi^\ast(Z(f)))
- kE_n^\prime$ is effective.  This contradicts \ref{t1}, 
however, 
since  $(\phi^\prime)^\ast(\pi^\ast Z(f)) = 
(\pi^{\prime})^\ast(\phi^\ast Z(f))$ 
and by hypothesis 
$$
\mult_{E_Y}(\phi^\ast Z(f)) <  \gamma \epsilon n
$$
and hence $(\pi^{\prime})^\ast(\phi^\ast Z(f))- (\pi^\prime)^\ast \gamma
\epsilon n \tilde{E}_Y $  is not effective.  

Choose a positive integer $b$ sufficiently large so that 
\begin{eqnarray}
T_X \ts bA \,\,\,\,\,\mbox{is generated by global sections}.
\next
\label{derivative}
\end{eqnarray}
Suppose $s \in H^0\left(X,n(D- \epsilon A)\right)$. We claim that
\next
\begin{eqnarray}
\ord_\eta(s) \geq \frac{n\epsilon}{b}.
\label{order}
\end{eqnarray}
If not then there exists a  local differential operator $D$ or order $t < 
\frac{n\epsilon}{b}$ such that $D(s)$ does not vanish at a general point 
of $Y$.  Using the
theory of \cite{eln}, there exists a differential section 
$$
s^\prime \in H^0\left(X,nD- n\epsilon A  + tbA\right)
$$
with $s^\prime(\eta) \neq 0$.  This contradicts the fact that $Y$ is in
the base locus of 
$$
H^0\left(X,nD - n\epsilon A + tbA\right)
$$
since $-n\epsilon + tb < 0$.  Thus we may take $\gamma = \frac{1}{b}$ and
this clearly does not depend upon $\epsilon$ or $n$, 
establishing Lemma \ref{l1}.

\nonmbthm
\setcounter{theorem}{6}
\medskip
\nmbthm
\init
\noindent
{\bf Proof of Lemma \ref{lemma2}}  Suppose $\eta \in Y$ is a general, hence
smooth, point
and let $\{x_1,\ldots,x_d\}$ be a system of local parameters at $\eta$ such
that $\ii_{Y} \subset \oo_X$ is generated on an open set
$U \supset \eta$ by $x_1, \ldots, x_c$.  Let $m_\eta \subset
\oo_{\eta,X}$ denote the maximal ideal of functions vanishing at $\eta$
and let
$$
U_{k,n} = \oo_{\eta,X}/m_\eta^{\delta kn}.
$$
We define $V_{k,n}(\epsilon) \subset U_{k,n}$ to
be the vector subspace generated by the images of the following collection
of monomials in the local parameters $\{x_1,\ldots,x_d\}$:
\next
\begin{eqnarray}
\left\{M(x_1,\ldots,x_c)M^\prime(x_{c+1},\ldots,x_d)| \,\, \deg(M) < 
\gamma \epsilon nk \,\,\,\mbox{and}\,\,\deg(M^\prime) < \delta nk\right\}.
\label{monomials}
\end{eqnarray}
By \ref{fifth}, we know that for $n$ sufficiently large the evaluation map
followed by the projection to the quotient
$$
\pi_{k,n}: H^0\left(X,kn(D + \epsilon A)\right) \ra U_{k,n}
$$
is surjective.  Choose a vector subspace $V^\prime_{k,n}(\epsilon)
\subset H^0\left(X,kn(D + \epsilon A)\right)$ such that
\next
\begin{eqnarray}
\pi_{k,n}: V^\prime_{k,n}(\epsilon)
 \ra V_{k,n}(\epsilon) \,\,\,\,\mbox{is an isomorphism}.
\label{iso}
\end{eqnarray}
Note that by \ref{fifth} we can assume that for $n$ sufficiently
large $n(D + \epsilon A)$ generates $\delta n$--jets at $\eta$ and
consequently we may assume that
\next
\begin{eqnarray}
V^{\prime}_{k,n}(\epsilon)
 \subset {\rm Sym}^k\left(H^0\left(X,n(D + \epsilon A)\right)\right).
\label{diso}
\end{eqnarray}

We claim that the restriction map
\next
\begin{eqnarray}
\phi_{k,n}: V^\prime_{k,n}(\epsilon) \ra
H^0\left(kY_n(\epsilon),
kn(D + \epsilon A)\right) \,\,\,\,\mbox{is injective}.
\label{inj}
\end{eqnarray}
By Lemma \ref{l1} it is sufficient to show that the restriction
$$
\psi_{k,n}: V^\prime_{k,n}(\epsilon) \ra
H^0\left(\gamma \epsilon knY,kn(D + \epsilon A)\right) 
$$
is injective.  If $\psi_{k,n}$ were not injective, there would exist
a non--zero polynomial 
$P(x_1,\ldots,x_d)$ in the monomials \ref{monomials} such that 
$P \in \ii_{\gamma \epsilon knY}$.  
Choose a monomial
$x^{a_1}\ldots x^{a_d}$ whose coefficient in $P$ is non--zero and such that
$a_1 + \ldots +a_c$ is minimal. Consider
$$
Q(x_1,\ldots,x_d) = 
\frac{\partial^{a_1}}{\partial x_1^{a_1}}\cdots
\frac{\partial^{a_c}}{\partial x_c^{a_c}}P(x_1,\ldots,x_d).
$$
By choice of $(a_1,\ldots,a_c)$, all terms in $Q(x_1,\ldots,x_d)$ vanish
identically along $Y$ except those coming from terms of type
\next
\begin{eqnarray}
\frac{\partial^{a_1}}{\partial x_1^{a_1}}\cdots
\frac{\partial^{a_c}}{\partial x_c^{a_c}}\left(x_1^{a_1}\ldots x_c^{a_c})
M^\prime(x_{c+1},\ldots,x_d\right).
\label{pst}
\end{eqnarray}
But no non--zero
sum of functions of the type appearing in \ref{pst} can vanish 
identically along $Y$ since the functions $x_{c+1},\ldots,x_{d}$ form a 
transcendence basis for the function field of $Y$.  Hence
\next
\begin{eqnarray}
\ord_Y(P(x_1,\ldots,x_d)) \leq a_1+ \ldots + a_c < \gamma \epsilon kn,
\label{ordery}
\end{eqnarray}
contradicting the assumption that $P \in \ii_{\gamma \epsilon kn Y}$.

We conclude from \ref{iso}  that 
$$
\dim\left(\phi_{k,n}\left(V_{k,n}^\prime(\epsilon)\right)\right) 
\geq \frac{(\gamma\epsilon nk)^{\dim Y}
(\delta nk)^{\codim (Y,X)}}{d!}
$$
By \ref{diso} each section $s \in V_{k,n}^\prime(\epsilon)$ vanishes along
$Z_n(-\epsilon)$ to order at least $k$ and
consequently vanishes along $kZ_n(-\epsilon)$.
Thus we can take $W_{k,n}(\epsilon) = 
\phi_{k,n}\left(V_{k,n}^\prime(\epsilon)\right)$ 
in Lemma \ref{lemma2}, with the
implied constant being $\frac{1}{d!}$, which does not depend upon $\epsilon$
as desired; the fact that no section of $W_{k,n}(\epsilon)$ 
vanishes along $Y$ to order
$\geq \gamma \epsilon kn$ is guaranteed by \ref{ordery}.

\nonmbthm
\setcounter{theorem}{7}
\nmbthm
\init
\medskip
\noindent
{\bf Proof of Corollary \ref{cor}}  As remarked above, in order for the
statement of Corollary \ref{cor} to make sense, we will assume that 
$D \in {\rm Div}(X) \ts \qq$ and $\epsilon \in \qq$.  We will also
write
$$
\rho = a\epsilon.
$$
  Assuming that $B$ is effective, and
for later purposes, very ample, consider the exact sequence
\next
\begin{eqnarray}
H^0\left(kY_n(\rho),kn(D + \rho A - a B)\right) & \ra&
H^0\left(kY_n(\rho),kn(D + \rho A)\right)  \nonumber \\
&\stackrel{\psi}{\ra}& 
H^0\left(aknB \cap kY_n(\rho),kn(D + \rho A)\right). 
\label{ex}
\end{eqnarray}
By hypothesis, we can write 
$$
akn B = \bigcup_{i=1}^{akn} B_i
$$
with each $B_i$ linearly equivalent to $B$.  
Let 
$H_W^0\left(B_i\cap kY_n(\rho),
kn(D + \rho A)\right) $ denote the 
image of $W_{k,n}(\rho)$, defined in Lemma \ref{lemma2}, in 
$H^0\left(B_i\cap kY_n(\rho),kn(D + \rho A)\right)$.  Since
$\dim (V_{k,n}(\rho)) = \dim({\rm ker}(\psi))$ 
we have, substituting into \ref{ex},
\next
\begin{eqnarray}
\dim(V_{k,n}(\rho)) \geq
\dim(W_{k,n}(\rho))- \sum_{i=1}^{akn}
h^0_W\left(B_i\cap kY_n(\rho),kn(D + \rho A)\right) .
\label{ext}
\end{eqnarray}
By Lemma \ref{lemma2}, it is sufficient, in order to establish Corollary
\ref{cor}, to prove that 
$$
\sum_{i=1}^{akn}
h_W^0\left(B_i\cap kY_n(\rho),
kn(D + \rho A)\right) \leq O\left(
a^{\codim (Y,X)+1}(nk)^{\dim X}\right).
$$
Thus for a general hypersurface $B$ we need to establish
\next
\begin{eqnarray}
h_W^0\left(B\cap kY_n(\rho),
kn(D + \rho A)\right) \leq O\left(
a^{\codim (Y,X)}(nk)^{\dim X-1}\right)
\label{ls}
\end{eqnarray}
with the implied constant depending only on $X$, $B$, and $\epsilon$.

We first eliminate the possibility the $\dim(Y) = 0$.  In this
case, $Y = P$ is a point, disjoint from the rest of the base locus 
$\base(D - \epsilon A)$.  
This is impossible, however, as \ref{thirdy} in combination with Lemma
\ref{l1} would establish that it is possible to lift a section
$$
s^\prime \in H^0(kY_n(\epsilon),kn(D-\epsilon A))
$$
with order of vanishing $< \gamma \epsilon kn$ at $P$, violating \ref{order}.

Assume now that $\dim Y \geq 1$.  
Let $X_{k,n}(\rho) = B \cap kY_n(\rho)$ and choose general
divisors $D_1,\ldots,D_c \in |B|$ where $c = \dim Y - 1$.
We have an exact sequence
\begin{eqnarray*}
H^0\left(X_{k,n}(\rho),kn(D+\rho A - aD_1))\right) 
&\stackrel{f_{k,n}}{\longrightarrow}& 
H^0\left(X_{k,n}(\rho),kn(D + \rho A)\right) \\
&\ra& H^0\left(aknD_1 
\cap X_{k,n}(\rho),kn(D+ \rho A)\right) .
\end{eqnarray*}
We claim that 
$$
{\rm Image}(f_{k,n}) \cap H_W^0\left(X_{k,n}(\rho),
kn(D+\rho A)\right) 
= \emptyset.
$$
Indeed, choosing the representative $D_1$ passing through a general point
$\eta \in Y$, we see that every section in ${\rm Image}(f_{k,n})$ has order
 of vanishing at least $akn$ along $\eta$. On the other hand, if 
$s \in H_W^0\left(X_{k,n}(\rho),kn(D+\rho A)\right)$ then
$\mult_\eta(s) < \gamma a \epsilon kn$; 
indeed, by \cite{Fu} Corollary 12.4, the 
multiplicity of a section at a point does not change under restriction to
a general hyperplane section.  We can, of course, assume that $\gamma \epsilon
 < 1$.
Hence, considering the exact sequence \ref{ex}, we obtain
$$
H_W^0\left(X_{k,n}(\rho),kn(D+\rho A)\right) \simeq
H_W^0\left(X_{k,n}(\rho) \cap  aknD_1,
kn(D + \rho A)\right),
$$
where the subscript $W$ again indicated sections coming from 
$W_{k,n}(\rho)$ by restriction.
It follows, using the same argument as after \ref{ext}, that 
$$
h_W^0\left(X_{k,n}(\rho),kn(D+ \rho A )\right) \leq 
akn h_W^0\left(X_{k,n}(\rho) \cap  D_1,kn(D+\rho A)\right).
$$
Proceding inductively gives
\next
\begin{eqnarray}
&h_W^0\left(X_{k,n}(\rho),
kn(D+\rho A)\right) \leq  & \nonumber \\
&(akn)^c h_W^0\left(X_{k,n}(\rho) \cap D_1 \cap \ldots \cap D_c,
kn(D+\rho A)\right).&
\label{is3}
\end{eqnarray}
Since $D_1,\ldots,D_c$ are general, the intersection $D_1  \cap  \ldots
\cap D_c \cap X_{k,n}(\rho)$ 
will contain only general points of $Y$ and thus
\ref{is3}, along with \cite{Fu} Corollary 12.4 and Example 4.3.4, implies
that
\next
\begin{eqnarray}
h_W^0\left(B \cap kY_n(\rho),kn(D+\rho A)\right) \leq 
(akn)^c \deg_B(Y_\red) \ell_Y(\oo_{kY_n(\rho)}).
\label{ut}
\end{eqnarray}

Thus in order to conclude the proof of Corollary \ref{cor}, we need to
control the length $\ell_Y(\oo_{kY_n(\rho)})$ as a function of $\rho = 
\epsilon a$.  We now reintroduce the rational parameter $a = p/q$,
and we will show that for any $0 < \frac{p}{q} < 1$ we
can find $n$ sufficiently large so that
\next
\begin{eqnarray}
\ell_Y\left(\oo_{kY_n\left(\frac{p\epsilon}{q}\right)}\right) 
\leq O\left(\left(\frac{p\epsilon kn}{q}\right)^{\codim(Y,X)}\right)
\label{yt}
\end{eqnarray}
where the implied constant depends only on $X$ and $\epsilon$.  
Combining \ref{ut} and \ref{yt} shows \ref{ls}, with appropriate divisibility
restrictions on $k$ and $n$ according to the rational number $p/q$,
and concludes the proof of Corollary \ref{cor}. 

\medskip

In order to prove \ref{yt}, 
we need to examine how the ideal sheaves $\ii_{kY_n(\epsilon)}$ 
depend upon $k$, $n$, and $\epsilon$
but only at a general point of $Y$ since we are only interested in
controlling the generic length of $\oo_{kY_n(\epsilon)}$.
Choose $\epsilon^\prime$ so that 
$0 < \epsilon^\prime - \epsilon \ll \epsilon$.
We begin by showing that for fixed $n$, we can find $\beta > 0$ such that
\next
\begin{eqnarray}
\ii_{\beta \epsilon^\prime knY} \subset \ii_{kY_n(\epsilon^\prime)}  
\,\,\,\,\mbox{generically} 
\,\,\,\mbox{for all}\,\,\,k \gg 0.
\label{asym}
\end{eqnarray}
In order to establish \ref{asym} we begin by fixing a $\beta > 0$ such that
$$
\ii_{\beta \epsilon^\prime nY} \subset \ii_{Y_n(\epsilon^\prime)},\,\,\,\,\,
\mbox{generically}.
$$
Such a $\beta$ must exist because, by assumption $\ii_{Y_n(\epsilon^\prime)}$
is generically $\ii_Y$--primary.  
We have 
$$
\ii_{kY_n(\epsilon^\prime)} \supset \ii_{Y_n(\epsilon^\prime)}^k 
\supset \ii_{\beta \epsilon^\prime
nY}^k \,\,\,\,\mbox{generically}\,\,\,\mbox{for all}\,\,\, k> 0.
$$
But, generically, $\ii_{\beta \epsilon^\prime nY}^k = 
\ii_{\beta \epsilon^\prime kn Y}$ for
all $k$ sufficiently large and thus we have established \ref{asym}.

Next, in order to allow $\epsilon^\prime$ to vary, we will need to let $n$
vary as well.  We will show the 
following: for any positive integer $r$ and for all $k > 0$
\next
\begin{eqnarray}
\ii_{kY_n(\epsilon_1)} \cdot \ldots \cdot \ii_{kY_n(\epsilon_r)} \subset
\ii_{kY_{rn}\left(\frac{\epsilon_1+\ldots+\epsilon_r}{r}\right)} \,\,\,\,
\mbox{generically};
\label{half}
\end{eqnarray}
here the $\epsilon_i$ are arbitrary rational numbers less than or equal
to $\epsilon^\prime$.  
Choose a birational map $\pi: V \ra X$ so that $V$ is smooth and so 
that $\pi$ resolves the base loci of the $r+1$ complete
linear series
$$
|n(D -\epsilon_1 A)|,\ldots, |n(D-\epsilon_r A)|, 
\left|rn\left(D  -
\frac{\epsilon_1+\ldots+\epsilon_r}{r} A\right)\right|.
$$
We let $E_1+F_1,\ldots,
E_r+F_r,$ and $E+F$ denote the associated exceptional divisors.  
The decomposition of the exceptional divisors is determined as follows:
we choose $E_1, \ldots, E_r, E$ to contain all exceptional components
whose  image under $\pi$ is contained in $Y$ and 
$F_1,\ldots, F_r, F$ are formed from all remaining components.  Thus,
no component of $F_1,\ldots,F_r,F$ is mapped into $Y$ by $\pi$.   
Then, generically,
 we have
\begin{eqnarray*}
\pi_\ast \oo_Y(-E_i) = \ii_{Y_n(\epsilon_i)} , \,\,\,\,1\leq i \leq r, \\
\pi_\ast \oo_Y(-E) = \ii_{Y_{rn}\left(\frac{\epsilon_1+\ldots+\epsilon_r}{r}
\right)},
\end{eqnarray*}
since, as we observed above with the scheme $Z_n(\epsilon)$, 
these ideal sheaves do not depend on the
resolution $\pi: V \ra X$.  

We claim that 
\next
\begin{eqnarray}
E_1 + \ldots + E_r \geq E.
\label{ine}
\end{eqnarray}
To see this, suppose $s_i \in H^0(X,n(D -\epsilon_i A))$ 
for $1 \leq i \leq r$.  It follows that
\next
\begin{eqnarray}
s_1 \ts \ldots \ts s_r \in  
H^0\left(X,rn\left(D -\frac{\epsilon_1+ \ldots +
\epsilon_r}{r} A\right)\right).
\label{hy}
\end{eqnarray}
Since the sections $s_i$ were arbitrary, 
lifting \ref{hy} to $V$ establishes that
$$
(E_1 + \ldots E_r) + (F_1 +\ldots F_r) \geq E + F:
$$
indeed \ref{hy} shows that the
sum of the common base loci of $|n(D - \epsilon_i A)|$ is at least
as large as the base locus of $\left|rn\left(D-
\frac{\epsilon_1+ \ldots + \epsilon_r}{r} A\right)\right|$.  
By definition of $E_i, F_i$, $E$, and $F$ this establishes \ref{ine}.
We can conclude that for all 
$k>0$
$$
\pi_\ast(\oo_V(-kE_1)) \cdot \ldots \cdot \pi_\ast(\oo_V(-kE_r)) \subset
\pi_\ast(\oo_V(-kE)), 
$$
and this is precisely \ref{half}, with the inclusion holding off of
$\supp(Y_{n,k}) \cap \supp(W_{n,k})$.

\medskip

In order to use \ref{half} to establish \ref{asym}, we apply \ref{half}
as follows.  Choose a rational number $\frac{p}{q} < 1$ and let 
$$
\epsilon_1 = \ldots = \epsilon_p = \epsilon^\prime
$$
and
$$
\epsilon_{p+1} = \ldots = \epsilon_q  = 
\frac{p\epsilon - p\epsilon^\prime}{q-p}.
$$
Applying \ref{half},
$$
\ii_{kY_n(\epsilon_1)} \cdot \ldots \cdot \ii_{kY_n(\epsilon_q)} \subset
\ii_{kY_{qn}\left(\frac{p\epsilon}{q}\right)}\,\,\,\,\,\mbox{generically for
all}\,\,\, k > 0.
$$
By the choice of $\epsilon_1, \ldots, \epsilon_q$ this gives
\next
\begin{eqnarray}
\ii_{kY_n(\epsilon^\prime)}^p \subset
\ii_{kY_{qn}\left(\frac{p\epsilon}{q}\right)} \,\,\,\,\,\mbox{generically
for all $k >0$}.
\label{ugh}
\end{eqnarray}
Using \ref{asym}, \ref{ugh} becomes
\next
\begin{eqnarray}
\ii_{\beta p\epsilon^\prime knY} \subset
 \ii_{kY_{qn}\left(\frac{p\epsilon}{q}\right)},\,\,\,\,\,\mbox{generically
for all $k \gg 0$}.
\label{hugh}
\end{eqnarray}
Now we will use \ref{hugh} in order to estimate lengths and establish \ref{yt},
concluding the proof of Corollary \ref{cor}. 
We have, by \ref{hugh} above,
\next
\begin{eqnarray}
\ell_Y\left(\oo_{kY_{qn}\left(\frac{p\epsilon}{q}\right)}\right) 
\leq \ell_Y(\oo_{\beta p \epsilon^\prime knY}).
\label{this}
\end{eqnarray}
Subsituting $n^\prime = qn$ into \ref{this} gives
\next
\begin{eqnarray}
\ell_Y\left(\oo_{kY_{n^\prime}\left(\frac{p\epsilon}{q}\right)}\right) 
\leq \ell_Y\left(
\oo_{\frac{\beta p \epsilon^\prime kn^\prime Y}{q}}\right).
\label{that}
\end{eqnarray}
A quick computation shows that 
$$
\ell_Y\left(\oo_{\frac{\beta p \epsilon^\prime knY}{q}}\right) = 
O\left(\left(\frac{\beta
\epsilon^\prime pkn^\prime}{q}\right)^{\codim (Y,X)}\right).
$$
This is precisely \ref{yt}, with $\epsilon^\prime$ in place of $\epsilon$
and with an additional factor of $\beta$ on the right hand side.
Since $\epsilon^\prime$ is arbitrarily close to $\epsilon$ 
and since $\beta$ is a fixed constant depending only on $\epsilon$ 
this establishes
Corollary \ref{cor}.

\medskip
\nonmbthm
\setcounter{theorem}{8}
\nmbthm
\init
\noindent
{\bf Proof of Lemma \ref{lemma3}}  
Choose $\epsilon_1, \epsilon_2 > \epsilon$ 
so that $D - \epsilon_1 A, D +\epsilon_2  A \in {\rm Div}(X) \ts \qq$.  
Choose
a birational map $\pi: V \ra X$, with $V$ smooth, resolving the base loci of 
$|n(D - \epsilon_1 A)|$ and $|n(D + \epsilon_2 A)|$
for appropriately divisible $n$.  Let $E_1$ and $E_2$ be the corresponding
exceptional divisors.  By assumption, $E_2$ contains no component
dominating $Y$ while $E_1$ contains at least one such component.  As in the 
proof of Corollary \ref{cor} we will write
$$
E_1 = E_1(Y) + F_1
$$
where every component of $E_1(Y)$ maps into $Y$ and no component of 
$F_2$ is mapped into $Y$ by $\pi$.  We claim that 
\next
\begin{eqnarray}
\pi_\ast \oo_V(-kF_1) = \ii_{kW_n(\epsilon_1)}, \,\,\,\mbox{for all}\,
\,  k > 0
\label{equal}
\end{eqnarray}
To prove \ref{equal} let $U = X \backslash Y$ and observe that 
$$
\pi_\ast \oo_V(-kF_1)|U = \ii_{kW_n(\epsilon_1)}|U.
$$
But since no component of $F_1$ is mapped into $Y$ it follows
that $\pi_\ast \oo_V(-kF_1)$ is the largest ideal sheaf whose restriction
to $U$ is equal to $\pi_\ast \oo_V(-kF_1)|U$; since 
$\ii_{kW_n(\epsilon_1)}$ was chosen to be maximal, having no embedded
components along $\supp(Y)$, this establishes \ref{equal}.


Choose a $\pi$--exceptional divisor $F \subset V$, containing no component
of $E_1(Y)$, 
such that
\next
\begin{eqnarray}
F > F_1
\label{choice}
\end{eqnarray}
Let $\jj = \pi_\ast \oo_V(-F)$.
Using \ref{equal} and \ref{choice}
establishes that
\next
\begin{eqnarray}
\jj^k \subset
\ii_{kW_n(\epsilon_1)},\,\,\,\,\forall k > 0.
\label{fst}
\end{eqnarray}
Moreover, if $r > 0$, we have by \ref{fst}
\next
\begin{eqnarray}
\jj^{kr} \subset \ii_{kW_n(\epsilon_1)}^r
\subset \ii_{kW_{rn}(\epsilon_1)}, \,\,\,\,\mbox{for all}\,\,\, r > 0
\label{ght}
\end{eqnarray}
and this is precisely Lemma \ref{lemma3} in the case where $m = 0$.

Fix an $m >0$ and choose a positive integer $r_m$
and a birational map
$$
\pi_m: X_m \ra X
$$
with $X_m$ smooth, resolving the base loci of the following complete linear
series:
$$
\left|(2^{m}-1)r_mn(D - \epsilon_1 A)\right|,
\left|n(D + \epsilon_2 A)\right|,
\left|2^m r_mn\left(D                    
-\epsilon_1 A + 
\frac{\epsilon_1 + \epsilon_2}{2^m} A\right)\right|.
$$
Let $G_1, G_2$, and $G^\prime$ 
denote the part of the exceptional divisor in the
three cases which does not lie over $Y$.
We claim that
\next
\begin{eqnarray}
G_1 + r_mG_2 \geq G^\prime.
\label{yes}
\end{eqnarray}
As we have seen above after \ref{hy}, 
this will be true provided that for any sections
$$
s_1 \in H^0(X, (2^m-1)r_mn(D-\epsilon_1A)), s_2 \in 
H^0(X, n(D +\epsilon_2A))
$$
we have
$$
s_1 \ts s_2^{\ts r_m} \in 
H^0\left(X, 2^m r_mn\left(D -\epsilon_1 A +
\frac{\epsilon_1 + \epsilon_2}{2^m} A\right)\right)
$$
which is an elementary computation of line bundles.  
Pushing \ref{yes}, multiplied by $k$, down to $X$ and using
\ref{ght}  yields
\next
\begin{eqnarray}
\jj^{kr_m} \cdot \ii_{kW_{(2^{m}-1)r_mn(-\epsilon_1)}} 
\subset \ii_{kW_{2^mr_mn}\left(-\epsilon_1+\frac{\epsilon_1
+\epsilon_2}{2^m}\right)}, \,\,\,\,\mbox{for all}\,\,\,
k \gg 0.
\label{rck}
\end{eqnarray}
We now replace $\epsilon_1$  with  $\beta$, 
which may necessitate increasing $n$ in 
order to clear denominators, and observe that 
$$
\ii_{kW_l(a)} \subset \ii_{kW_l(b)}, \,\,\,\,\forall k,l
$$
whenever $a > b$.  Thus we conclude from \ref{rck} that for any 
$\beta < \epsilon$
$$
\jj^{kr_m} \cdot \ii_{kW_{(2^{m}-1)r_mn(-\beta)}} 
\subset \ii_{kW_{2^mr_mn}\left(-\beta+\frac{\epsilon}{2^m}\right)}, 
\,\,\,\,\mbox{for all}\,\,\,
k \gg 0,
$$
provided, as observed above, that $r_m$ is sufficiently divisible.
This concludes the proof of Lemma \ref{lemma3} and hence the proof of 
Theorem \ref{seshlimit}.  If $D \in {\rm Div}(X) \ts \qq$ 
then the same argument
works for $\beta = 0$.

\section{A theorem of Fujita}

In this section we turn to a theorem of Fujita regarding a numerical
form of Zariski decomposition.  The techniques of this paper are so close
to those of \cite{DEL} that it is only natural to apply them in this setting.
If $L$ is a big line bundle, we let $v(L)$ denotes its 
volume as defined in \cite{DEL}, so $v(L)$ measures the asymptotic growth of
$h^0(X,nL)$.  

\begin{thm}[Fujita]  
Suppose $X$ is a smooth projective variety and $L$ a big line
bundle on $X$.  For any $\epsilon > 0$ there exists a birational map
$\pi: Y \ra X$ and {\bf Q}--divisors $A$ and $E$ such that
$$
\pi^\ast L \equiv A + E
$$
where $A$ is ample and $v(A) \geq v(L) - \epsilon$.  
\label{fuj}
\end{thm}

We will attack this problem using the language of moving Seshadri constants
developed in \S 1.  Consider the following inductive construction.  Choose
a point 
$$
x_1 \in X \backslash \base(L)
$$
such that $\epsilon_m(x_1,L)$ is maximal. 
Next choose $\delta_1 \ll \epsilon_m(x_1,L)$ so that
$$
\epsilon_1 = \epsilon_m(x_1,L) - \delta_1 \in \qq.
$$
 Let $\pi_1: X_1 \ra X$ be the blow--up 
of $x_1$ with exceptional divisor $E_1$ and let
$$
L_1 = \pi_1^\ast L(-\epsilon_1E_1).
$$
Next we choose a point $x_2$ where $\epsilon_m(x_2,L_1)$
is maximal and let 
$\delta_2 \ll \epsilon_m(x_2,L_1)$ so that $\epsilon_2 = 
\epsilon_m(x_2,L_1) - \delta_2 \in \qq$.  Let
$\pi_2: X_2 \ra X$ be the blow up of $X$ at both $x_1$ and
$x_2$ and let
$$
L_2 = \pi_2^\ast L(-\epsilon_1 E_1  - \epsilon_2 E_2).
$$
This procedure can be repeated inductively to yield  pairs $(X_n,L_n)$ for
all positive $n$.  

\begin{lm} Suppose $\eta$ is a general point of $X$, which we also view
as a point of $X_n$.  Then
$$
\lim_{n \rightarrow \infty} \epsilon_m(\eta,L_n) = 0.
$$
\label{aty}
\end{lm}

\noindent
{\bf Proof of Lemma \ref{aty}}  If $\lim_{n \rightarrow \infty} 
\epsilon_m(\eta,L_n) \neq 0$ this means that there exists $\alpha > 0$
such that for $k$ suitably large and divisible
$kL$ generates $k\alpha$--jets at arbitrarily many points $x_1,x_2,
\ldots$, assuming that $\delta_i$ is sufficiently small relative to 
$\epsilon_i$.  This is clearly impossible as $v(L)$ bounds the number
of such points.

\begin{lm} Let
 $D$ be a very ample line bundle on $X$ and
assume that $L_n(-\alpha \pi_n^\ast D)$ is not in the effective cone.  
Then
$$
h^0(X_n,kL_n) \leq \alpha O(k^{\dim X}),
$$
with the implied constant depending only on $D$ and $L$.  
\label{bty}
\end{lm}

\nmbthm
\init
\addtocounter{nmb}{1}

\noindent
{\bf Proof of Lemma \ref{bty}}  
Consider the exact sequence
\begin{eqnarray*}
0 \ra H^0(X_n, kL_n(-\alpha \pi_n^\ast D)) \ra H^0(X_n,kL_n) \ra
H^0(k\alpha \pi_n^\ast D,kL_n)
\label{hit}
\end{eqnarray*}
where $k$ is sufficiently large and divisible and where $k\alpha \pi_n^\ast D$
is represented by the pull-back of $k\alpha$ general
smooth hypersurfaces linearly
equivalent to $D$.  Then we obtain
\begin{eqnarray*}
h^0(X_n,kL_n) &\leq & k\alpha  h^0(\pi_n^\ast D,kL_n) \nonumber \\
              & = &   k\alpha  h^0(D,kL) \nonumber \\
              & = &  \alpha O(k^{\dim X}).
\end{eqnarray*}
In the last formula, the implied constant clearly depends only on $D$ and $L$.


\medskip

\noindent
{\bf Proof of Theorem \ref{fuj}}
Choose $\alpha > 0$.  By Lemma \ref{aty} 
we can find $n$ so that $L_n - \alpha \pi_n^\ast D$
is not in the  effective cone.
Fix $k_0$ so that $|k_0L|$ generates $k_0\epsilon_i$ jets at $x_i$
for $1 \leq i \leq n$; we assume here that $k_0$ is sufficiently divisible
so that $k_0 \epsilon_i \in \zz$ for all $i$.  Thus the sequence
$$
0 \ra H^0\left(X, mk_0L \ts \cap_{i=1}^n m_{x_i}^{mk_0\epsilon_i}\right) 
\ra H^0(X,mk_0L) \ra
H^0\left(mk_0L \ts \oo_X/\cap_{i=1}^n m_{x_i}^{mk_0\epsilon_i}\right)  \ra 0
$$
is exact for all positive integers $m$; here $m_{x_i}$ is the 
maximal ideal sheaf of the point $x_i$.  Also by construction and 
Lemma \ref{bty}
\begin{eqnarray*}
h^0\left(X, mk_0L \ts \cap_{i=1}^n m_{x_i}^{mk_0\epsilon_i}\right) 
& = & h^0(X_n, mk_0 L_n) \\
& \leq &  \alpha c (mk_0)^{\dim X}
\end{eqnarray*}
with $c$ depending only on  $D$ and $L$.  
Let $\pi: Y \ra X$  be a birational map such that $Y$ is smooth and
$$
\pi^\ast k_0L \sim M + F
$$
with $M$ base point free.  Then  $v(M/k_0) \geq v(L) - \alpha c$.  Since
$\alpha$ was arbitrary, this establishes Theorem \ref{fuj}: note that here
$M/k_0$ is only big and nef rather than ample but this can be changed with
an arbitrarily small perturbation.

\medskip

Finally, we would like to pose a difficult and important question which
 we will phrase as a conjecture.  According
to Theorem \ref{seshlimit} if $D$ is effective and $A$ ample then
$D- \alpha A$ reaches the boundary of the effective cone exactly when
$$
\epsilon_m(\eta,D-\alpha A) = 0
$$
for a general point $\eta$.  Our question then is the following:
\begin{conjecture}  Suppose $D \in {\rm Div}(X) \ts \qq$ is a divisor class 
on the boundary of the effective cone.  Then $\epsilon_\eta(D) = 0$ for 
a general point $\eta \in X$.
\label{conj}
\end{conjecture}
Usint the results of this paper, one can establish under the hypotheses of
Conjecture \ref{conj} that
$$
\epsilon_m(\eta,D) = 0
$$
but it is not clear whether or not $\epsilon_m(\eta,D) = \epsilon(\eta,D)$.
If the answer is yes, then one has a numerical criterion to determine the 
effective cone.  

Conjecture \ref{conj} essentially boils down to the following 
concrete, computational question.
Suppose $L$ is an effective line bundle and suppose $\pi_n:Y_n \ra X$
resolves the base locus of $|nL|$ with $Y_n$ is smooth.
Writing $(\pi_n)^\ast nL = M_n + E_n$
is
$$
\lim_{n \ra \infty} \frac{M_n^{\dim X -1} \cdot E_n}{n^{\dim X}} = 0?
$$
In other words, does the pair $(M_n,E_n)$ approximate a Zariski decomposition
not only in the sense of Theorem \ref{fuj} but also in a stronger sense of
becoming closer and closer to being perpindicular?

\begin{tabbing}
Department of Mathematics and Statistics\\
University of New Mexico\\
Albuquerque, New Mexico 87131\\
{\em Electronic mail:} nakamaye@math.unm.edu
\end{tabbing}

\end{document}